\documentclass[12pt]{article}

\usepackage{amsmath, epsfig, cite}
\usepackage{amssymb}
\usepackage{amsfonts}
\usepackage{latexsym}
\usepackage{graphicx}
\usepackage{color}
\usepackage{ifpdf}
\usepackage{CJK}

\newtheorem{thm}{Theorem}[section]
\newtheorem{prop}[thm]{Proposition}

\newtheorem{defi}[thm]{Definition}

\newtheorem{problem}[thm]{Problem}

\numberwithin{equation}{section}

\makeatletter \@addtoreset{equation}{section} \makeatother

\begin{document}
\rule{0cm}{1cm}

\begin{center}
{\Large\bf Rainbow connection number of dense graphs }
\end{center}

\begin{center}
{\small Xueliang Li$^1$, Mengmeng Liu$^1$, Ingo Schiermeyer$^2$\\
$^1$Center for Combinatorics and LPMC-TJKLC\\
Nankai University, Tianjin 300071, China\\
lxl@nankai.edu.cn, liumm05@163.com\\
$^2$Institut f\"{u}t Diskrete Mathematik und Algebra\\
Technische Universit\"{a}t Bergakademie Frieiberg\\
09596 Freiberg, Germany\\
Ingo.Schiermeyer@tu-freiberg.de }
\end{center}

\begin{center}
\begin{minipage}{120mm}
\begin{center}
{\bf Abstract}
\end{center}

{\small An edge-colored graph $G$ is rainbow connected, if any two
vertices are connected by a path whose edges have distinct colors.
The rainbow connection number of a connected graph $G$, denoted
$rc(G)$, is the smallest number of colors that are needed in order
to make $G$ rainbow connected. In this paper we show that $rc(G)\leq
3$, if $|E(G)|\geq {{n-2}\choose 2}+2$, and $rc(G)\leq 4$, if
$|E(G)|\geq {{n-3}\choose 2}+3$. These bounds are sharp.\\[2mm]
{\bf AMS classification 2010}: 05C15, 05C40. }
\end{minipage}
\end{center}

\section{Introduction}

We use \cite{BM} for terminology and notation not defined here and
consider finite and simple graphs only.

An edge-colored graph $G$ is called rainbow connected, if any two
vertices are connected by a path whose edges have different colors.
The concept of rainbow connection in graphs was introduced by
Chartrand et al. in \cite{C}. The {\it rainbow connection number} of
a connected graph $G$, denoted $rc(G)$, is the smallest number of
colors that are needed in order to make $G$ rainbow connected. The
rainbow connection number has been studied for several graph
classes. These results are presented in a recent survey \cite{Li}.
Rainbow connection has an interesting application for the secure
transfer of classified information between agencies (cf. \cite{E})

In \cite{KS} the following problem was suggested:
\begin{problem}
For every $k,1\leq k \leq n-1,$ compute and minimize the function
$f(n,k)$ with the following property: If $|E(G)|\geq f(n,k)$, then
$rc(G)\leq k.$
\end{problem}

The authors of \cite{KS} got the following results:
\begin{prop}
$f(n,k)\geq {{n-k+1}\choose 2}+(k-1).$
\end{prop}

\begin{prop}
$f(n,1)={n \choose 2}, f(n,n-1)=n-1, f(n,n-2)=n.$
\end{prop}

\begin{thm}
Let $G$ be a connected graph of order $n$ and size $m$. If ${{n-1}
\choose 2}+1\leq m \leq {n \choose 2}-1$, then $rc(G)=2.$
\end{thm}

Hence $f(n,2)={{n-1} \choose 2}+1$. In this paper we will show that
$f(n,3)={{n-2} \choose 2}+2$ and $f(n,4)={{n-3} \choose 2}+3$.

\section{Main results}

\begin{defi}
Let $G$ be a connected graph. The distance between two vertices $u$
and $v$ in $G$, denoted by $d(u,v)$, is the length of a shortest
path between them in $G$. Distance between a vertex $v$ and a set $S
\in V(G)$ is defined as $d(v,S)=min_{x\in S}d(v,x)$. The $k$-step
open neighborhood of a set $S \in V(G)$ is defined as $N^k(S)=\{x
\in V(G)|d(x,S)=k \}, k \in \{0,1,2\cdots\}.$ When $k=1,$ we may
omit the qualifier ``$1$-step" in the above name and the superscript
$1$ in the notation. The neighborhood of a vertex $v$ in
$\overline{G}$, denoted by $\overline{N}(v)$, is defined as
$\overline{N}(v)=\{x|xv \notin E(G)\}$.
\end{defi}

We will present a new and shorter proof to the following result in
\cite{KS}.
\begin{thm}\label{thm1}
Let $G$ be a connected graph of order $n \geq 3$. If ${{n-1} \choose
2} +1 \leq |E(G)| \leq {n \choose 2}-1$, then $rc(G)= 2.$
\end{thm}

\begin{pf}
Our proof will be by induction on $n$. For $n=3$, we have
$f(n,n-1)=n-1=2={{3-1} \choose 2}+1$, for $n=4$, we have
$f(n,n-2)=n=4={{4-1} \choose 2}+1.$ So we may assume $n\geq 5.$

Since $|E(G)| \leq {n \choose 2} -1,$ we have $1 \leq \delta(G) \leq
n-2.$ Choose a vertex $w \in V(G)$ with $ d(w) = \delta(G)$ and set
$d(w)=n-2-t$ with $0 \leq t \leq n-3.$ Let $H = G -w,$ then $|E(H)|
\geq {{n-1}\choose 2}+1-d(w)={{n-2}\choose
2}+n-2+1-(n-2-t)={{n-2}\choose 2}+1+t={({n-1)-1}\choose
2}+1+t\geq{{n-2}\choose 2}+1.$ Hence, $H$ is connected, otherwise
$E(H)<{{n-2} \choose 2} +1$.

Now let $\overline{N}(w)=\{v_1,v_2,\cdots,v_t,v_{t+1}\}$.

Claim: $N(v_i)\bigcap N(w) \neq \emptyset$ for $1\leq i\leq t+1.$

\begin{pf}
Suppose $N(v_i)\bigcap N(w) = \emptyset$  for some $i,1\leq i\leq
t+1$. Then $d(v_i)\leq (t+1)-1=t$, thus
$E(\overline{G})=|\overline{N}_{H}(v_i)|+|\overline{N}_{G}(w)|\geq(n-t-2)+(t+1)=n-1>n-2$,
a contradiction, since $E(\overline{G})\leq n-2.$
\end{pf}

Hence for every vertex $v_i,$ there is a vertex $u_i \in N(w)$ such
that $u_iv_i \in E(G)$ for $1\leq i \leq t+1$. Let $H'$ be a
subgraph of $H$ with $V(H')=V(H)$ and
$E(H')=E(H)-\{v_1u_1,\cdots,v_tu_t\}$. Then $|E(H')|\geq {{n-2}
\choose 2}+1+t-t={{n-2} \choose 2} +1={{(n-1)-1} \choose 2} +1$. So
$H'$ is connected, by induction, $rc(H')\leq 2$. Now take a
$2$-rainbow colouring of $H'$. Let $c(v_{t+1}u_{t+1})=1$, then set
$c(v_iu_i)=1$ for $1\leq i\leq t$ and $c(e)=2$ for all edges $e$
which are incident with $w$. Then $G$ is $2$-rainbow connected.
\end{pf}
\begin{qed}
\end{qed}

The following are new results of this paper.
\begin{thm}\label{thm2}
Let $G$ be a connected graph of order $n \geq 4$. If $|E(G)| \geq
{{n-2} \choose 2} + 2$, then $rc(G)\leq 3.$
\end{thm}

\begin{pf}
Our proof will be by induction on $n$. For $n=4$, we have
$f(n,n-1)=n-1=3={{4-2} \choose 2}+2$, for $n=5$, we have
$f(n,n-2)=n=5={{5-2} \choose 2}+2.$ So we may assume $n\geq 6.$

By Theorem \ref{thm1}, we have $rc(G)\leq 2$ for $|E(G)|\geq {{n-1}
\choose 2}+1$. Hence we may assume $|E(G)|\leq {{n-1} \choose 2}$.
This implies $\delta(G)\leq
\frac{(n-1)(n-2)}{n}=n-3+\frac{2}{n}<n-2$.

Claim $1$: $diam(G)\leq 3.$

\begin{pf}
Suppose $diam(G)\geq 4$ and consider a diameter path
$v_1,v_2,\cdots,v_{D+1}$ with $D\geq 4$. Then $d(v_1)+d(v_4)\leq
n-2$ and $d(v_2)+d(v_5)\leq n-2$ implying $|E(G)|\leq {n \choose
2}-2(2n-3-(n-2))={n \choose 2}-2(n-1)={{n-2} \choose 2}-1<{{n-2}
\choose 2}+2$, a contradiction.
\end{pf}

Claim $2$: If $\delta(G)=1$, then $rc(G)\leq 3.$

\begin{pf}
Let $w$ be a vertex with $d(w)=\delta(G)=1$, let $H=G-w$. Then
$|E(H)|\geq {{n-2}\choose 2}+2-1={{n-2}\choose
2}+1={{(n-1)-1}\choose 2}+1$, hence $rc(H)\leq 2$ by Theorem
\ref{thm1}. Take a $2$-rainbow colouring for $H$ and set $c(e)=3$
for the edge incident with $w$. Then $rc(G)\leq 3.$
\end{pf}

Hence we may assume $\delta(G)\geq 2$. Let $w_1,w_2 \in V(G)$ with
$w_1w_2 \notin E(G)$. Suppose $N(w_1)\bigcap N(w_2)= \emptyset$. Let
$H=G-\{w_1,w_2\}$. Then $|E(H)|\geq {{n-2} \choose 2}+2-(n-2)={{n-3}
\choose 2}+1={{(n-2)-1} \choose 2}+1$. Thus $H$ is connected. Hence
$rc(H)\leq 2$ by Theorem \ref{thm1}. Consider a $2$-rainbow
colouring of $H$ with colours $1$ and $2$. Since $diam(G)\leq 3$,
there is a $w_1w_2$-path $w_1u_1u_2w_2$. Let $c(u_1u_2)=1$, then set
$c(w_1u_1)=2, c(w_2u_2)=3$ and $c(e)=3$ for all other edges incident
with $w_1$ or $w_2$. Then $G$ is $3$-rainbow connected.

Hence we may assume if $w_1,w_2 \in V(G),w_1w_2 \notin E(G)$, then
$N(w_1)\bigcap N(w_2)\neq \emptyset$. Choose a vertex $w$ with
$d(w)=\delta(G)$ and set $d(w)=n-2-t$ with $1\leq t \leq n-4$. As in
the proof of Theorem \ref{thm1}, there exist vertices $u_i \in N(w)$
such that $u_iv_i \in E(G)$ for $1\leq i\leq t+1$. Let $H=G-w$ and
$H'$ be a subgraph of $H$ with $V(H')=V(H)$ and
$E(H')=E(H)-\{u_1v_1,\cdots,u_{t-1}v_{t-1}\}$. Then $|E(H')|\geq
{{n-2} \choose 2}+2-(n-2-t)-(t-1)={{n-2} \choose 2}-n+5={{(n-1)-2}
\choose 2}+2$.

Hence, if $H'$ is connected, by induction, $H'$ is $3$-rainbow
connected. Now take a $3$-rainbow colouring of $H'$. Let
$c(u_iv_i)\in \{1,2\}$ for $i=t,t+1$. Then set $c(u_iv_i)=1$ for
$1\leq i\leq t-1$ and $c(e)=3$ for all edges $e$ incident with $w$.
Then $G$ is $3$-rainbow connected.

Claim $3$: If $H'$ is disconnected, then $H'$ has at most $2$
components and one of them is a single vertex.

\begin{pf}
Suppose, by way of contradiction, that $H'$ has $k\geq 3$
components, let $n_i$ be the number of vertices of the $i$th
component, thus $n_1+\cdots+n_k=n-1$, then

\begin{eqnarray*}
|E(H')| & \leq &
{n_1\choose 2}+\cdots+{n_k \choose 2}\\
&=& \displaystyle\sum_{i=1}^k\frac{n_i^2-n_i}{2}=\frac{1}{2}\left(\displaystyle\sum_{i=1}^kn_i^2-(n-1)\right)\\
&\leq & \frac{1}{2}\left[1+1+(n-1-2)^2-n+1\right]\\
&=& \frac{1}{2}\left(n^2-7n+12\right)\\
&<& {{n-3}\choose 2}+2,
\end{eqnarray*}
a contradiction. So $H$ has two component, that is $n_1+n_2=n-1$, if
$n_1 \geq 2,$ then
\begin{eqnarray*}
|E(H')| & \leq &
{n_1\choose 2}+{n_2 \choose 2}\\
&=&\frac{n_1^2+n_2^2-(n-1)}{2}\\
&\leq & \frac{1}{2}\left[2^2+(n-3)^2-n+1\right]\\
&=& \frac{1}{2}\left(n^2-7n+14\right)\\
&<& {{n-3}\choose 2}+2.
\end{eqnarray*}
thus complete the proof.
\end{pf}

Let $H_1=\{v\},H_2$ be two components of $H'$, we know that $v \in
N(w)$(otherwise $\delta(G)=1$), let
$N(v)=\{w,v_1,\cdots,v_{d_x-1}\}$. Obviously $d_x \leq t,$ and all
edges $v_iv, 1\leq i\leq d_x-1$ are deleting edges. Since
$|E(H_2)|=|E(H')|\geq {{(n-2)-1} \choose 2}+2$, $H_2$ is $2$-rainbow
connected by Theorem \ref{thm1}. Consider a $2$-rainbow colouring of
$H_2$ with colours $1,2$.  Set $c(vv_i)=3, 1\leq i\leq d_x-1,
c(wv)=1,c(e)=3$ for all edges incident with $w, c(e)=2$ for all
other deleting edges. Then for every $x \in V(G)\backslash{w}$,
there is a rainbow path between $w$ and $x$, and for every $x\in
N(w)$, there is a rainbow path $vwx$, for every $x\in
N^2(w)\backslash N(v)$, we know $xv \notin E(G)$, then $N(x)\bigcap
N(v)\neq \emptyset,$ that means there exist some $v_i, 1\leq i\leq
d_x-1$ with $v_i\in N(x)\bigcap N(v)$, i.e. there is a rainbow path
between $v$ and $x$. So $G$ is $3$-rainbow connected.
\end{pf}
\begin{qed}
\end{qed}

\begin{thm}\label{thm3}
Let $G$ be a connected graph of order $n \geq 5$. If $|E(G)| \geq
{{n-3} \choose 2} + 3$, then $rc(G)\leq 4.$
\end{thm}

\begin{pf}
We apply the proof idea from the proof of Theorem \ref{thm2}.

Our proof will be by induction on $n$. For $n=5$, we have
$f(n,n-1)=n-1=4={{5-3} \choose 2}+3$ and for $n=6$, we have
$f(n,n-2)=n=5={{6-3} \choose 2}+3$. So we may assume $n\geq 7$.

By Theorem \ref{thm2}, we have $rc(G)\leq 3$ for $|E(G)|\geq {{n-2}
\choose 2}+2$. Hence we may assume $|E(G)|\leq {{n-2} \choose 2}+1$.
This implies $\delta(G)\leq
\frac{(n-2)(n-3)}{n}=n-5+\frac{8}{n}<n-3$.

Claim $4$: $diam(G)\leq 4.$

\begin{pf}
Suppose $diam(G)\geq 5$ and consider a diameter path
$v_1,v_2,\cdots,v_{D+1}$ with $D\geq 5$. Then $d(v_i)+d(v_{i+3})\leq
n-2$ for $i=1,2,3$ implying $|E(G)|\leq {n \choose
2}-3(2n-3-(n-2))={n \choose 2}-3(n-1)={{n-3} \choose 2}-3<{{n-3}
\choose 2}+3$, a contradiction.
\end{pf}

Claim $5$: If $\delta(G)=1$, then $rc(G)\leq 4.$

\begin{pf}
Let $w$ be a vertex with $d(w)=\delta(G)=1$, let $H=G-w$. Then
$|E(H)|\geq {{n-3}\choose 2}+3-1={{n-3}\choose
2}+2={{(n-1)-2}\choose 2}+2$, hence $rc(H)\leq 3$ by Theorem
\ref{thm2}. Take a $3$-rainbow colouring for $H$ and set $c(e)=4$
for the edge incident with $w$. Then $rc(G)\leq 4.$
\end{pf}

Hence we may assume $\delta(G)\geq 2$.

Case $1:$ If there are $w_1,w_2 \in V(G), w_1w_2 \notin E(G)$, with
$N(w_1)\bigcap N(w_2)\neq \emptyset$ such that  $d(w_1)+d(w_2)\leq
n-3.$

Let $H=G-\{w_1,w_2\}$. Then $|E(H)|\geq {{n-3} \choose
2}+3-(n-4)={{n-4} \choose 2}+2+1={{(n-2)-2} \choose 2}+2$. We claim
that $H$ is connected. Or by the proof of Theorem \ref{thm2}, we
know $H$ has at most $2$ components and one of them is a single
vertex. Thus $\delta(G)=1$, a contradiction. Hence $rc(H)\leq 3$ by
Theorem \ref{thm2}. Consider a $3$-rainbow colouring of $H$ with
colours $1,2,3$. If there is a rainbow path $P=xyz$ of length $2$
between $N(w_1)$ and $N(w_2)$, where $x \in N(w_1),z \in N(w_2)$,
let $c(xy)=1,c(yz)=2$, then set $c(w_1x)=3, c(w_2z)=4$ and $c(e)=4$
for all other edges incident with $w_1$ or $w_2$. Then $G$ is
$4$-rainbow connected. If all paths of length $2$ between $N(w_1)$
and $N(w_2)$ are not rainbow connected, then we choose a path
$P=xyz$, where $x \in N(w_1),z \in N(w_2)$, let $c(xy)=c(yz)=1$,
then we keep the colour of all edges in $E(H)$ except $yz$. Then set
$c(yz)=4,c(w_1x)=2,c(w_2z)=3$ and $c(e)=4$ for all other edges
incident with $w_1$ or $w_2$. It's only need to check $G$ is
$4$-rainbow connected. Since $\delta \geq 2$, then there exists $v$
such that $c(w_1v)=4$. For every $w \in V(G)/N(w_1)$, there is a
rainbow path $P$ not including $yz$, Otherwise there is a rainbow
path of length $2$ between $N(w_1)$ and $N(w_2)$, so $w_1vPw$ is a
rainbow path. It's similar for $w_2$.

Case $2:$ For all $w_1, w_2 \in V(G)$ with $w_1w_2 \notin E(G),
N(w_1)\bigcap N(w_2)\neq \emptyset$.

We know that in this case $diam(G)=2$. Choose a vertex $w$ with
$d(w)=\delta(G)$ and set $d(w)=n-2-t$ with $2\leq t \leq n-4$.

As in the proof of Theorem \ref{thm1}, there exist vertices $u_i \in
N(w)$ such that $u_iv_i \in E(G)$ for $1\leq i\leq t+1$. Let $H=G-w$
and $H'$ be a subgraph of $H$ with $V(H')=V(H)$ and
$E(H')=E(H)-\{u_1v_1,\cdots,u_{t-2}v_{t-2}\}$. Then $|E(H')|\geq
{{n-3} \choose 2}+3-(n-2-t)-(t-2)={{n-3} \choose 2}-n+7={{(n-1)-3}
\choose 2}+3$.

If $H'$ is connected, by induction, $H'$ is $4$-rainbow connected.
Now take a $4$-rainbow colouring of $H'$. Let $c(u_iv_i)\in
\{1,2,3\}$ for $i=t-1,t,t+1$. Then set $c(u_iv_i)=1$ for $1\leq
i\leq t-2$ and $c(e)=4$ for all edges $e$ incident with $w$. Then
$G$ is $4$-rainbow connected.

If $H'$ is disconnected, we claim that $H'$ has at most $3$
components. Otherwise $|E(H')|<{{n-4}\choose 2} +3.$ If $H'$ has
exactly $3$ components $H_1,H_2,H_3,$ we may assume that $|H_3|\geq
|H_2|\geq |H_1|\geq 1, |H_1|+ |H_2|+ |H_3|=n-1$. If $|H_2|\geq 2,$
then $|E(H')|\leq {|H_1| \choose 2}+{|H_2| \choose 2}+{|H_3| \choose
2}\leq 1+{{n-4}\choose 2}<{{n-4}\choose 2} +3.$ So $|H_1|=|H_2|=1,$
let $V(H_1)=\{u_1\}, V(H_2)=\{u_2\}$ and $u_1, u_2 \in N(w).$ Then
$|E(H_3)|\geq {{n-4} \choose 2}+3 \geq {{(n-3)-1} \choose 2}+3.$
Hence, by Theorem \ref{thm1}, $H_3$ is $2$-rainbow connected. Now
consider a $2$-rainbow colouring of $H'$ with $1, 2.$ Set
$c(wu_1)=1, c(wu_2)=2, c(e)=4$ for all edges incident with $w,
c(e)=3$ for all edges incident with $u_1, u_2 $ except $wu_1,wu_2$.
$c(e)=1$ for all other deleting edges. Then $G$ is $4$-rainbow
connected.

If $H'$ has exactly $2$ components $H_1,H_2$. we may assume that
$|H_2|\geq |H_1|\geq 1$. First, $|H_1|=1,$ let $V(H_1)=\{u_1\}$ and
$u_1 \in N(w)$. Then $|E(H_2)|\geq {{n-4} \choose 2}+3 \geq
{{(n-2)-2} \choose 2}+3.$ Hence, by Theorem \ref{thm2}, $H_2$ is
$3$-rainbow connected. Now consider a $3$-rainbow colouring of $H'$
with $1,2,3.$ Set $c(wu_1)=1, c(e)=4$ for all edges incident with
$w$ or $u_1$ except $wu_1, c(e)=2$ for all other deleting edges.
Then $G$ is $4$-rainbow connected. Second, $|H_1|\geq 2$. Since
$n\geq 7$, then $|H_2|\geq 3$. Thus we have
\begin{eqnarray*}
|E(H_1)| & \geq &
{{n-4}\choose 2}+3-{|H_2| \choose 2}\\
&= & \frac{1}{2}\left[|H_1|^2-3|H_1|+4\right]+|H_1||H_2|-3|H_2|-2|H_1|+7\\
&\geq&{{|H_1|-1 }\choose 2}+1+3(n-4)-3(n-1)+|H_1|+7\\
&\geq & {{|H_1|-1 }\choose 2}+1
\end{eqnarray*}

Similar $|E(H_2)| \geq {{|H_2|-1 }\choose 2}+1$. Hence both
$H_1,H_2$ are $2$-rainbow connected. Consider a $2$-rainbow
colouring of $H'$ with $1, 2.$ Set $c(wv)=4$, for all $v \in V(H_1),
c(wv)=3$, for all $v \in V(H_2), c(uv)=4$ for all $u \in
V(H_2)\bigcap N(w), v \in V(H_1)\bigcap N^2(w), c(uv)=3$ for all $u
\in V(H_1)\bigcap N(w), v \in V(H_2)\bigcap N^2(w), c(e)=1$ for all
other edges. Then $G$ is $4$-rainbow connected.

Case $3:$ For all $w_1,w_2 \in V(G), w_1w_2 \notin E(G)$, such that
$d(w_1)+d(w_2)\geq n-2.$

We know that in this case $diam(G)=3$. Choose a vertex $w$ with
$d(w)=\delta(G)$ and set $d(w)=n-2-t$ with $2\leq t \leq n-4$.

Subcase $2.1:$ $N^3(w)=\emptyset.$

It is same as the proof of Case $2$.

Subcase $2.2:$ $N^3(w)\neq \emptyset.$

For every $u \in N^3(w), wu \notin E(G)$ and $N(w)\bigcap
N(u)=\emptyset$, then $d(w)+d(u)= n-2,$ that is $N(u)=N^2(w)\bigcup
N^3(w) \backslash\{u\}$. Let $N(w)=\{u_1,\cdots,u_{n-t-2}\},
N^2(w)=\{v_1,\cdots,v_p\},p\geq1,N^3(w)=\{v_{p+1},\cdots,v_{t+1}\}$.

If $p=1$, $v_1$ is a cut vertex and $G[N^2(w)\bigcup N^3(w)]$ is a
complete graph. Let $H_1,H_2$ be two blocks of $G$, we may assume
that $H_2$ is a complete graph. We know that $H_1$ is a bipartite
graph $K_{2,n-t-2}$, then $H_1$ is $4$-rainbow connected. Now give a
$4$-coloring of $H_1$ as follows:
$$
c(e)=\left \{
\begin{array}{ll}
j+1,& \mbox {if $e=u_iw$ if $i\in \{3j+1,3j+2,3j+3\}$ for $0\leq j\leq 2$},\\
 4, & \mbox {if $e=u_iw$ for $i>9$},\\
 \mbox{$i$ mol $3$}, & \mbox {if $e=v_1u_i$ for $i\leq 9$},\\
 3,& \mbox{if $e=v_1v_i$ for $i>9$}.
\end{array}
\right.
$$
and let $c(e)=4$ for $e \in E(H_2)$. Then $G$ is $4$-rainbow
connected.

If $p=2$, let $H_1=G[w\bigcup N(w)\bigcup
N^2(w)],H_2=G[N^2(w)\bigcup N^3(w)]$, then $|H_1|+|H_2|=n+2,
|H_1|\geq 5, |H_2|\geq 3$, then
\begin{eqnarray*}
|E(H_1)| & \geq &
{{n-3}\choose 2}+3-{|H_2| \choose 2}\\
&= & \frac{1}{2}\left[|H_1|^2-5|H_1|+6\right]+2+|H_1||H_2|-3|H_1|-5|H_2|+13\\
&\geq &{{|H_1|-2 }\choose 2}+2+5(n+2-5)-5(n+2)+2|H_1|+13\\
&\geq & {{|H_1|-2 }\choose 2}+2
\end{eqnarray*}
if $|H_1|\geq 6$. Then $H_1$ is $3$-rainbow connected. Consider a
$3$-rainbow colouring of $H_1$ with $2,3,4.$  Set $c(e)=1$ for all
$e \in E(H_2)$, then $G$ is $4$-rainbow connected. When $|H_1|=5$,
set $c(wu_1)=4, c(wu_2)=3, c(u_1v_1)=2, c(u_2v_2)=1, c(e)=1$ for all
$e \in E(H_2)$, then $G$ is $4$-rainbow connected.

Now we may assume that $p\geq 3$. For every $v_i \in N^2(w)$, there
is a vertex $u_i \in N(w)$ such that $u_iv_i \in E(G).$ Let $H$ be
the graph be deleting $w$ and edges $u_iv_i$ for $v_i \in
N^2(w)\backslash\{v_1,v_2,v_3\}$ and edges $v_1v_i$ for $p+1\leq
i\leq t+1$, then $|E(H)|\geq {{n-4 \choose 2}}+3.$ If $H$ is
connected, by induction $H$ is $4$-rainbow connected. Consider a
$4$-rainbow colouring of $H$ with $1,2,3,4.$ Let $c(u_iv_i)\in
\{1,2,3\}$ for $i=1,2,3$. We may assume that $c(u_1v_1)=1$. Set
$c(e)=4$ for all edges incident with $w$, $c(v_1v_i)=2, p+1\leq
i\leq t+1$, $c(v_iu_i)=3$ for all other edges between $N(w)$ and
$N^2(w)$. Then $G$ is $4$-rainbow connected.

If $G$ is disconnected, as proof in case $2$, $H$ has at most $3$
components. If $H$ has exactly $3$ components with two single
vertices in $N(w)$, denoted by $u_1,u_2$, then $H_1$ is $2$-rainbow
connected. Consider a $2$-rainbow colouring of $H_1$ with $1,2.$ Let
$c(uv_1)=1$ for $u \in N(w)$, set $c(wu_1)=c(wu_2)=3,c(e)=4$ for all
edges incident with $w$, $c(v_1v_i)=2, p+1\leq i\leq t+1$, $c(e)=3$
for all edges incident with $u_1$ except $wu_1$, $c(e)=4$ for all
edges incident with $u_2$ except $wu_2, c(e)=1$ for all the
remaining edges. Then $G$ is $4$-rainbow connected.

If $G$ has exactly two components $H_1,H_2$. First $|H_1|=1$, let
$V(H_1)=\{u_1\} \subseteq N(w)$, $H_2$ is $3$-rainbow connected by
the proof of case $2$. Now consider a $3$-rainbow colouring of $H_2$
with $1,2,3.$ Let $c(uv_1)=1$ for $u \in N(w)\backslash \{u_1\}$,
set $c(wu_1)=3,c(e)=4$ for all edges incident with $w$,
$c(v_1v_i)=2, p+1\leq i\leq t+1$, $c(e)=2$ for all edges incident
with $u_1$ except $wu_1, c(e)=1$ for all the remaining edges. Then
$G$ is $4$-rainbow connected. Second $|H_1|\geq 2$, $H_1,H_2$ both
$2$-rainbow connected by the proof of case $2$. We may assume that
$H_2$ contains $N^3(w)$. Now consider a $2$-rainbow colouring of
$H_1,H_2$ with $1,2.$ Set $c(wv)=3$ for all $v \in V(H_1)$,$c(wv)=4$
for all $v \in V(H_2)$ $c(v_1v_i)=4, p+1\leq i\leq t+1, c(e)=3$ for
all the remaining edges. Then $G$ is $4$-rainbow connected.
\end{pf}
\begin{qed}
\end{qed}

\end{document}